\documentclass[11pt]{article}
\usepackage{fullpage, amsmath, amssymb, amsthm, epsfig, color, graphicx}

\def\eqbd{\mathop{{:}{=}}}
\def\bdeq{\mathop{{=}{:}}}
\def\Q{\mathbb{Q}}
\def\R{\mathbb{R}}
\def\P{\mathbb{P}}

\def\C{\mathbb{C}}

\def\A{{\bf{A}}}
\def\a{{\bf{a}}}

\def\qed{\hfill {$\Box$}}

\def\two_n{{\underbrace{2{\times}2{\times}\cdots{\times}2}_{n \; terms}}}

\newenvironment{@abssec}[1]{%
     \if@twocolumn
       \section*{#1}%
     \else
       \vspace{.05in}\footnotesize
       \parindent .2in
         {\bfseries #1. }\ignorespaces
     \fi}
     {\if@twocolumn\else\par\vspace{.1in}\fi}
\newenvironment{keywords}{\begin{@abssec}{Key words}}{\end{@abssec}}
\newenvironment{AMS}{\begin{@abssec}{AMS subject classification}}{\end{@abssec}}

\newtheorem{theorem}{Theorem}

\newtheorem{remark}[theorem]{Remark}
\newtheorem{lemma}[theorem]{Lemma}
\newtheorem{corollary}[theorem]{Corollary}
\newtheorem{conjecture}[theorem]{Conjecture}
\newtheorem{problem}[theorem]{Problem}

\newtheorem{example}[theorem]{Example}
\newcommand{\ignore}[1]{}

\title{Hyperdeterminantal relations among \\ symmetric principal minors}
\author{ Olga Holtz \\ Mathematics Department \\ University of California \\
Berkeley, CA 94720 \and
Bernd Sturmfels  \\ Mathematics Department \\ University of California \\
Berkeley, CA 94720}

\begin{document}
\maketitle

\begin{abstract}
The principal minors of a symmetric $n{\times}n$-matrix
form a vector of length $2^n$. We characterize these vectors
in terms of algebraic equations derived from the
$ 2{\times}2{\times}2$-hyperdeterminant.
\end{abstract}



\section{Introduction}

The principal minors of a real symmetric $n{\times}n$-matrix
$A$ form a vector of length $2^n$.
This vector is denoted $A_*$, and its entries are
indexed by subsets $I$ of $[n]\eqbd \{1,2,\ldots,n\}$. Namely, $A_{I}$ denotes the
minor of $A$ whose rows and columns are
indexed by $I$.  This includes the $0 {\times} 0$-minor
 $\,A_{\emptyset} = 1$.
The aim of this paper is to give an
 algebraic characterization of all vectors
 in $\R^{2^n}$ which arise in this form.
  Our question can be rephrased as follows.
We write $\R^{\binom{n+1}{2}}$ for the space of 
real symmetric $n{\times}n$-matrices.
 Our aim is to determine the image of the
{\em principal minor map} $\,\phi : \R^{\binom{n+1}{2}} \rightarrow \R^{2^n},\, A \mapsto A_* $.
  If $n = 2$ then the image of $\phi$
 is characterized by the trivial equation $A_\emptyset = 1$ and
 the one inequality
 \begin{equation}
 \label{HF1} A_{\emptyset} \cdot A_{12} \quad \leq \quad A_{1} \cdot A_{2}.
 \end{equation}
 
 For $n \geq 2$,  the entries of the symmetric matrix $A = (a_{ij})$ are
 determined up to sign by their principal minors
 of size $1{\times}1$ and $2{\times}2$,
 in view of the relation
$\, a_{ij}^2 = A_{i} A_{j} - A_{ij} A_\emptyset$.

\begin{remark}\label{closure}
The image of the principal minor map $\phi$  is a closed subset
of $\R^{2^n}$.
The same holds for the map
$\,\phi_\C : \C^{\binom{n+1}{2}} \rightarrow \C^{2^n}\,$
which takes a complex symmetric matrix to its  principal minors.
\end{remark}

\begin{proof}
Let $A_*$  be a vector in $\R^{2^n}$ (resp.~in $\C^{2^n}$) that
lies in the closure of the image of $\phi$ (resp.~the image of $\phi_\C$).
Consider any sequence $\,\{A^{(k)}\}_{k \geq 0}\,$ of symmetric $n{\times}n$-matrices whose principal
minors tend to $A_*$ as $k\to \infty$. This sequence must be uniformly
bounded, since the diagonal entries of $A^{(k)}$ have prescribed limits, and so do the
magnitudes of all off-diagonal entries. By compactness, we can extract a convergent
subsequence $A^{(k_j)}$, whose limit $A$ therefore has principal minors $A_*$.
\end{proof}

 Remark \ref{closure} implies that the image of $\,\phi_C\,$
is a complex algebraic variety in $\C^{2^n}$. That variety
 is the object of study in this paper. Equivalently, we seek to
 determine all polynomial equations that are valid
 on the image of the map $\phi : \R^{\binom{n+1}{2}} \rightarrow \R^{2^n}$.

  For $n = 3$, the matrix $A$ has six distinct entries, so the
 seven non-trivial minors $A_{I}$ must satisfy one polynomial equation.
 Expanding the determinant $A_{123}$,
 the desired equation is found to be
$$
\left( A_{123}-A_{12} A_{3} -A_{13 } A_{2}-A_{23} A_{1} +2A_{1 }
A_{2} A_{3}  \right)^2 \quad = \quad
4 \cdot (A_{1} A_{2} - A_{12})(A_{2} A_{3} - A_{23})(A_{1}A_{3}- A_{13}),
$$
Our point of departure is the observation that
this equation coincides with the
{\em hyperdeterminant} of format $2 {\times} 2 {\times} 2$.
Namely, after homogenization using
$A_{\emptyset} = 1$, this equation is equivalent to
 \begin{equation}
\begin{array}{ll}
 &  A_{\emptyset}^2 A_{123}^2  \,+\, 
 A_{1}^2 A_{23}^2 \,+\,   A_{2}^2 A_{13}^2 \,+\, A_{3}^2 A_{12}^2
\, + \,\,4 \cdot A_{\emptyset} A_{12} A_{13} A_{23} \,\,+\,\,
 4 \cdot A_{1} A_{2} A_{3} A_{123}
\\ &   -\,\,  2 \cdot  A_{\emptyset}   A_{1} A_{23} A_{123}
 \,-\,  2 \cdot  A_{\emptyset}   A_{2} A_{13} A_{123}
 \,-\,  2 \cdot A_{\emptyset}   A_{3} A_{12} A_{123}
 \,-\,  2 \cdot A_{1} A_{2} A_{13} A_{23} \\ &
 \,-\,  2 \cdot A_{1} A_{3} A_{12} A_{23}
\, - \,2 \cdot  A_{2} A_{3} A_{12} A_{13}
 \quad = \quad 0 .  \end{array}   \label{hyperrel}
 \end{equation}
  Replacing subsets of $\{1,2,3\}$ by binary strings
 in $\{0,1\}^3$, and thus setting $\,A_{\emptyset} = a_{000},
A_{1} = a_{100}$,
$ \ldots, A_{123} = a_{111}$, we see that this polynomial has
 the full symmetry group of the $3$-cube, and it does indeed coincide with the
 familiar formula for the hyperdeterminant
(see \cite[Prop.~14.B.1.7]{GKZbook}).
The {\em hyperdeterminantal relations} on $A_{*}$
(of format $2 {\times} 2 {\times} 2 $)
are obtained from~(\ref{hyperrel}) by substituting
\begin{equation}
\begin{array}{ll}
& A_{\emptyset} \mapsto A_{I} \,\, , \quad
A_{123} \mapsto A_{I \cup \{j_1,j_2,j_3\}}  \,,
\\ &
A_{1} \mapsto A_{I \cup \{j_1\}}\, ,\,\,
A_{2} \mapsto A_{I \cup \{j_2\}} \,,\,\,
A_{3} \mapsto A_{I \cup \{j_3\}}  \\ &
A_{12} \mapsto A_{I \cup \{j_1,j_2\}}\, ,\,\,
A_{13} \mapsto A_{I \cup \{j_1,j_3\}} \,,\,\,
A_{23} \mapsto A_{I \cup \{j_2,j_3\}} 
\end{array}  \label{substitute}
\end{equation}
for any subset $I \subset [n]$ and
any elements $j_1,j_2,j_3 \in [n] \backslash I $.
Our first result states:
\begin{theorem}
\label{first}
The principal minors of a symmetric matrix satisfy the
hyperdeterminantal relations.
\end{theorem}

The proof of this theorem is presented in Section~\ref{sec_background}, after
reviewing some relevant background from matrix theory and the history of our
problem. We also comment on connections to probability theory.
A key question is whether the converse to Theorem~\ref{first} holds,
i.e., whether every vector of length $2^n$ which satisfies the
hyperdeterminantal relations arises from the principal minors
of a (possibly complex) symmetric $n {\times} n$-matrix.
The answer to this question is ``not quite but almost''.
A useful counterexample to keep in mind is the following
  point on the hyperdeterminantal locus.

\begin{example}
\label{caution} \rm
For  $n\geq 4$, define the vector $\,A_{*} \in \R^{2^n}\,$ by
 $A_{\emptyset} = 1$, $\, A_{123\cdots n} = -1$ and
$A_{I} = 0$ for all other subsets $I $ of $[n]$. Then  $A_{*}$ satisfies
all hyperdeterminantal relations. It also satisfies all the
{\em Hadamard-Fischer inequalities,} which are the generalization of the inequality
(\ref{HF1}) to arbitrary $n$:
\begin{equation}
 \label{HF2} A_{I \cap J} \cdot A_{I \cup J}
\quad \leq \quad A_{I} \cdot A_{J} \quad
 \qquad \hbox{for all} \quad I,J \subseteq [n].
 \end{equation}
 These inequalities hold for positive semidefinite symmetric matrices
 (see, e.g.,~\cite{EngelSchneider}).
 
 Nonetheless, the vector $A_*$ is not the vector of principal minors of any
symmetric $n {\times} n$-matrix. To see this, we note that a symmetric
$n{\times}n$-matrix $A$ which has all diagonal entries $A_i$ and all
principal $2{\times}2$-minors $A_{ij}$ zero must be the zero matrix, so its
determinant $A_{123\cdots n}$ would be zero.
\end{example}

We spell out different versions of a converse to Theorem~\ref{first}
in the later sections. In Section~\ref{sec_converse1} we derive a
converse under a genericity  hypothesis,
and in Section~\ref{sec_bighyperdet} we discuss larger
hyperdeterminants and derive a converse using
 so-called condensation relations. The ultimate converse
would be an explicit list of generators
for the prime ideal $P_n$ of the algebraic variety ${\rm image}(\phi_\C)$.
In Section~\ref{sec_landsberg} we conjecture that
the ideal $P_n$ is generated by quartics,
namely, the orbit of the hyperdeterminantal relations
under a natural group action. Section~\ref{sec_ideal}
verifies this conjecture for $n=4$.
Theorem~\ref{first} together with these converses resolves the
{\em Symmetric Principal Minor Assignment Problem}
which was stated as an open question in
Problem 3.4 of~\cite{HoltzSchneider} and in Section 3.2 of~\cite{Wagner}.

\section{Matrix Theory and Probability}  \label{sec_background}

This work is motivated by a number of recent results and problems from matrix theory and
probability. Information about the principal minors of a given matrix is crucial in many
matrix-theoretic settings. Of interest may be their exact value, their sign, or inequalities
they satisfy. Among these problems are detection of $P$-matrices~\cite{TsatLi} and of
$GKK$-matrices~\cite{gkktau, HoltzSchneider}, counting spanning trees of a graph~\cite{GorniTZ}
and the inverse multiplicative eigenvalue problem~\cite{Friedland2}. The
{\em Principal Minor Assignment Problem,} as formulated in~\cite{HoltzSchneider}, is to determine whether a given vector $A_*$ of length $2^n$ is realizable as the vector of all principal minors
of some $n {\times} n$-matrix $A$. 
Very recently, Griffin and Tsatsomeros  gave an algorithmic solution to this
problem~\cite{GriffinTsat1, GriffinTsat2}. Their work gives an algorithm, which,  under a certain ``genericity'' condition, either outputs a solution matrix or determines that none exists.
Our approach offers a more conceptual algebraic solution in the case of
symmetric matrices.

In probability theory, information about principal minors is important in {\em determinantal
point processes.} Determinantal processes arise naturally in several fields, including
quantum mechanics of fermions~\cite{Macchi}, eigenvalues of random matrices, random spanning
trees and nonintersecting paths (see~\cite{Houghetal} and references therein).
A determinantal  point process on a locally compact measure space $(\Lambda,\mu)$ is determined
by a kernel  $K(x,y)$ so that the joint intensities of the process can be written as
$\det(K(x_i,x_j))$.  In particular, if $\Lambda$ is a finite set and $\mu$ is the counting measure
on $\Lambda$, then $K$ reduces to a $|\Lambda| {\times} |\Lambda|$-matrix, whose principal minors give
the joint densities of the process. The matrix $K$ is not necessarily Hermitian (or real symmetric),
even though very often it is.
 
 In the theory of negatively correlated random variables~\cite{Pemantle, Wagner}, principal
minors of a real $n {\times} n$-matrix give  values of a function $\omega : 2^{[n]} \to [0,\infty)$
on the Boolean algebra  $2^{[n]}$ of all subsets of  $[n] = \{1,2,\ldots,n\}$. The function
$\omega$ must be non-negative and must meet the following negative correlation condition. Suppose
$y_1,\ldots,y_n$ are indeterminates, and consider the {\em partition function}
$$ Z(\omega; y)\,\, \eqbd \,\,\sum_{I\subseteq [n]} \omega(I) \cdot y^I, \qquad
\hbox{where} \quad y^I \, \eqbd \, \prod_{i\in I} y_i. $$
For any positive vector $y = (y_1,\ldots,y_n)$, this determines a probability measure
$\mu=\mu_y$ on $2^{[n]}$ by
$$ \mu(I) \,\, \eqbd \,\,{\omega(I) y^I \over Z(\omega;y) } \qquad \hbox{\rm for all} \;\;
 I \subseteq [n]. $$
The {\em atomic random variables} of this theory are given by $X_i(I)\eqbd 1$ if $i\in I$
and $0$ otherwise. Their expectations and covariances are
$$ \langle X \rangle \,\eqbd\, \sum_{I\subseteq [n]} X(I) \mu(I) \qquad
\mbox{and} \qquad {\rm Cov}(X,Y)\, \eqbd\,
\langle XY \rangle -  \langle X \rangle \langle Y \rangle, $$
 and the {\em negative correlation hypothesis} requires that
$$ {\rm Cov} (X_i,X_j)\leq 0  \qquad \hbox{\rm for all}\;\; y>0, \;\; i\neq j.  $$
Wagner \cite{Wagner} asks how to characterize all functions $\omega$ satisfying these conditions
and arising from some matrix $A$, i.e., such that $\omega(I)=A_I$ is the minor of $A$
with columns and rows indexed by the subset $I \subseteq [n]$.
This application to probability theory is one of the motivations for our
algebraic approach to the principal minor
assignment problem, namely, the characterization of algebraic relations among principal minors. 
In this paper we restrict ourselves to the symmetric case.

We now recall a basic fact about Schur complements
 (e.g. from~\cite{BrualdiSchneider}). The {\em Schur complement}
of an invertible principal  submatrix $H$ in a matrix $\,A\,$ is the
 matrix  $\,AH\eqbd E-FH^{-1}G\,$ where
 $$ A \,\,\bdeq \,\,\left( \begin{array}{cc} E & F \\ G & H  \end{array} \right) .$$  
 The Schur complement is the result
of Gaussian elimination applied to reduce the submatrix $F$ to zero using the rows of $H$.
The principal minors of the Schur  complement satisfy {\em Schur's identity}
\begin{equation} \qquad (AH)_\alpha\,\,\,= \,\,\, {A_{I \cup  \alpha} \over A_I}
\quad \qquad \hbox{for all subsets} \,\, \alpha \,\subseteq \, [n] \backslash I,
  \label{schur_id}
\end{equation}
assuming that $H$ is the principal submatrix of $A$ with rows and columns indexed by  $I$.

The hyperdeterminantal relations of format $2 {\times} 2 {\times} 2$ are now
 derived from Schur's identity (\ref{schur_id}):

\smallskip

\noindent {\em Proof of Theorem \ref{first}. \/ }
The validity of the relation (\ref{hyperrel}) for symmetric
$3 {\times} 3$-matrices is an easy direct calculation.
Next suppose that $A$ is a symmetric $n {\times} n$-matrix
all of whose principal minors $A_I$ are non-zero.
The hyperdeterminantal relation for
$I \cup \{j_1,j_2,j_3\}$ coincides with the
relation $(\ref{hyperrel})$ for the
principal $3 {\times} 3$-minor indexed by
$\{j_1,j_2,j_3\}$ in the Schur complement $AH$,
after multiplying by $A_I^4$ to clear denominators.
Here we are using Schur's identity (\ref{schur_id})
for any non-empty subset $\alpha$ of $\{j_1,j_2,j_3\}$.
Now, if $A$ is any symmetric matrix that has vanishing
principal minors then we write $A$ as the limit of a sequence
of matrices whose principal minors are non-zero. The
hyperdeterminantal relations hold for every matrix in the
sequence, and hence they hold for $A$ as well.
\qed

\vskip .2cm

In Section~\ref{sec_bighyperdet} we offer an alternative proof of Theorem~\ref{first},
by showing that the vector $A_*$ satisfies the hyperdeterminantal
relations  of higher-dimensional formats $2 {\times} 2 {\times} \cdots {\times} 2$.
At this point we note that the hyperdeterminantal relations
do not suffice even if all principal minors are non-zero.

\begin{example} \label{eighthyper}
\rm
For $n=4$ there are $8$ hyperdeterminantal relations,
one for each facet of the $4$-cube:
\begin{eqnarray*}
 &  A_{\emptyset}^2 A_{123}^2  \,+\,   A_{1}^2 A_{23}^2 \,+\,   A_{2}^2 A_{13}^2
 \,+\, A_{3}^2 A_{12}^2 \, + \,\,4  A_{\emptyset} A_{12} A_{13} A_{23} \,\,+\,\,
   4A_{1} A_{2} A_{3} A_{123} \,  -\,\,  2   A_{\emptyset}   A_{1} A_{23} A_{123}  \\ &
 \,-\,   2A_{\emptyset}   A_{2} A_{13} A_{123}   \,-\,   2A_{\emptyset}   A_{3} A_{12} A_{123}
\,-\,  2 A_{1} A_{2} A_{13} A_{23}  \,-\,  2 A_{1} A_{3} A_{12} A_{23} \, - \,2 A_{2} A_{3} A_{12} A_{13}
\,=\, 0,
\\ &  A_{\emptyset}^2 A_{124}^2  \,+\,   A_{1}^2 A_{24}^2 \,+\,   A_{2}^2 A_{14}^2
 \,+\, A_{4}^2 A_{12}^2 \, + \,\,4  A_{\emptyset} A_{12} A_{14} A_{24} \,\,+\,\,
   4A_{1} A_{2} A_{4} A_{124} \,  -\,\,  2   A_{\emptyset}   A_{1} A_{24} A_{124}  \\ &
 \,-\,   2A_{\emptyset}   A_{2} A_{14} A_{124}   \,-\,   2A_{\emptyset}   A_{4} A_{12} A_{124}
\,-\,  2 A_{1} A_{2} A_{14} A_{24}  \,-\,  2 A_{1} A_{4} A_{12} A_{24} \, - \,2 A_{2} A_{4} A_{12} A_{14}
\,=\,0, \\ &  A_{\emptyset}^2 A_{134}^2 \,+\,   A_{1}^2 A_{34}^2 \,+\,   A_{3}^3 A_{14}^2
 \,+\, A_{4}^3 A_{13}^2 \, + \,\,4  A_{\emptyset} A_{13} A_{14} A_{34} \,\,+\,\,
   4A_{1} A_{3} A_{4} A_{134} \,  -\,\,  2   A_{\emptyset}   A_{1} A_{34} A_{134}  \\ &
 \,-\,   2A_{\emptyset}   A_{3} A_{14} A_{134}   \,-\,   2A_{\emptyset}   A_{4} A_{13} A_{134}
\,-\,  2 A_{1} A_{3} A_{14} A_{34}  \,-\,  2 A_{1} A_{4} A_{13} A_{34} \, - \,2 A_{3} A_{4} A_{13} A_{14}
\,=\,0, \\ &  A_{\emptyset}^2 A_{234}^2  \,+\,   A_{2}^2 A_{34}^2 \,+\,   A_{3}^2 A_{24}^2
 \,+\, A_{4}^2 A_{23}^2 \, + \,\,4  A_{\emptyset} A_{23} A_{24} A_{34} \,\,+\,\,
   4A_{2} A_{3} A_{4} A_{234} \,  -\,\,  2   A_{\emptyset}   A_{2} A_{34} A_{234}  \\ &
 \,-\,   2A_{\emptyset}   A_{3} A_{24} A_{234}   \,-\,   2A_{\emptyset}   A_{4} A_{23} A_{234}
\,-\,  2 A_{2} A_{3} A_{24} A_{34}  \,-\,  2 A_{2} A_{4} A_{23} A_{34} \, - \,2 A_{3} A_{4} A_{23} A_{24}
\,=\,0,
\end{eqnarray*}
\begin{eqnarray*}
 &   A_{4}^2 A_{1234}^2  +   A_{14}^2 A_{234}^2 +   A_{24}^2 A_{134}^2
 + A_{34}^2 A_{124}^2  + 4  A_{4} A_{124} A_{134} A_{234} +
   4A_{14} A_{24} A_{34} A_{1234} \\ &  -  2   A_{4}   A_{14} A_{234} A_{1234} 
 -   2A_{4}   A_{24} A_{134} A_{1234}   -   2A_{4}   A_{34} A_{124} A_{1234} \\ & \qquad
-  2 A_{14} A_{24} A_{134} A_{234}  -  2 A_{14} A_{34} A_{124} A_{234}  - 2 A_{24} A_{34} A_{124} A_{134} \,\, = \,\, 0
, \\ &  A_{3}^2 A_{1234}^2  +   A_{13}^2 A_{234}^2 +   A_{23}^2 A_{134}^2
 + A_{34}^2 A_{123}^2  + 4  A_{3} A_{123} A_{134} A_{234} +
   4A_{13} A_{23} A_{34} A_{1234}   \\
& -  2   A_{3}   A_{13} A_{234} A_{1234} 
 -   2A_{3}   A_{23} A_{134} A_{1234}   -   2A_{3}   A_{34} A_{123} A_{1234}  \\ & \qquad
-  2 A_{13} A_{23} A_{134} A_{234}  -  2 A_{13} A_{34} A_{123} A_{234}  - 2 A_{23} A_{34} A_{123} A_{134} \,\,=\,\,0 ,\\
&  A_{2}^2 A_{1234}^2  +   A_{12}^2 A_{234}^2 +   A_{23}^2 A_{124}^2
 + A_{24}^2 A_{123}^2  + 4  A_{2} A_{123} A_{124} A_{234} +
   4A_{12} A_{23} A_{24} A_{1234} \\ &   -  2  A_{2}   A_{12} A_{234} A_{1234} 
 -   2 A_{2}   A_{23} A_{124} A_{1234}   -   2 A_{2}   A_{24} A_{123} A_{1234}  \\ & \qquad
-  2 A_{12} A_{23} A_{124} A_{234}  -  2 A_{12} A_{24} A_{123} A_{234} 
- 2 A_{23} A_{24} A_{123} A_{124} \,\,=\,\, 0,
 \\ &  A_{1}^2 A_{1234}^2  +   A_{12}^3 A_{134}^2 +   A_{13}^3 A_{124}^2
 + A_{14}^3 A_{123}^3  + 4  A_{1} A_{123} A_{124} A_{134} +
   4A_{12} A_{13} A_{14} A_{1234}  \\ & -  2   A_{1}   A_{12} A_{134} A_{1234} 
 -   2A_{1}   A_{13} A_{124} A_{1234}   -   2A_{1}   A_{14} A_{123} A_{1234}  \\ & \qquad
-  2 A_{12} A_{13} A_{124} A_{134}  -  2 A_{12} A_{14} A_{123} A_{134} 
- 2 A_{13} A_{14} A_{123} A_{124} \,\, = \,\, 0.
\end{eqnarray*}
The set of solutions to these equations has the correct codimension (five)
but it is too big. To illustrate this phenomenon, consider the
case when all minors of a given size have the same value:
$$
A_\emptyset \,:=\, x_0 \, , \,\,\,
 A_i \,:= \,x_1\,, \,\,\,
 A_{ij} \,:= \,x_2\,, \,\,\,
 A_{ijk} \,:= \,x_3\,, \,\,\,
 A_{1234} \,:=\, x_4 .$$
Under this specialization, the eight hyperdeterminants above reduce to
a system of two equations:
\begin{equation}
\label{comp1}
x_0^2 x_3^2-6 x_0 x_1 x_2 x_3+4 x_0 x_2^3+4 x_1^3 x_3-3 x_1^2 x_2^2 \,\,=\,\,
x_1^2 x_4^2-6 x_1 x_2 x_3 x_4+4 x_1 x_3^3+4 x_2^3 x_4-3 x_2^2 x_3^2 \,\,=
\,\, 0.
\end{equation}
The solution set to these equations in $\P^4$ is the union of
two irreducible surfaces, of degree ten and six respectively.
The degree ten surface is
extraneous and is gotten by requiring additionally that
\begin{equation}
\label{comp2}
 x_0 x_3^2 \,\, = \,\, x_1^2 x_4 .
\end{equation}
The degree six surface is our desired locus of principal minors.
It is defined by the two equations
\begin{equation}
\label{comp3}
     3 x_2^2-4 x_1 x_3+x_0 x_4 \,\, = \,\,
     4 x_2^3-6 x_1 x_2 x_3+x_0 x_3^2+x_1^2 x_4 \,\, = \,\, 0 .
\end{equation}
For a concrete numerical example
let us consider the symmetric $4 {\times} 4$-matrix
$$ A \quad = \quad \begin{pmatrix}
2 & 1 & 1 & 1 \\
1 & 2 & 1 & 1 \\
1 & 1 & 2 & 1 \\
1 & 1 & 1 & 2
\end{pmatrix} .$$
Its principal minors are
$\,(x_0,x_1,x_2,x_3,x_4) =  (1,2,3,4,5) \,$ and these
satisfy (\ref{comp3}) but not (\ref{comp2}).
On the other hand, the vector
$\, (x_0,x_1,x_2,x_3,x_4) =  (1,2,3,4,4) \,$
satisfies (\ref{comp1}) and (\ref{comp2}) but
not (\ref{comp3}). The corresponding vector
$\,A_* \in \R^{16}\,$ has all its entries non-zero and
satisfies the eight hyperdeterminantal relations but it
does not come from the principal minors of any
symmetric $4 {\times} 4$-matrix.
\qed
\end{example}

\section{A First Converse to Theorem~\ref{first}}   \label{sec_converse1} 

We now derive two additional classes of relations that the principal minors
of a symmetric $n {\times}n$ matrix $A =(a_{ij})$ must satisfy.
Throughout this section we set $A_\emptyset = 1$.
For $n\geq 4$ and for distinct $i$, $j$,
$k$, $l$, we can write the product
$\,8a^2_{ij}a^2_{ik}a^2_{il}a_{jk}a_{jl}a_{kl}\,$ in two different ways:
\begin{eqnarray}
&& \!\!\!\! \!(A_{ijk} - A_{ij} A_{k} - A_{ik } A_{j} - A_{jk} A_{i} + 2A_{i}  A_{j}  A_{k})
 (A_{ijl} - A_{ij}  A_{l} - A_{il}  A_{j} - A_{jl}  A_{i} + 2 A_{i}  A_{j}  A_{l}) \nonumber \\
&&  \times  (A_{ikl} - A_{ik}  A_{l} - A_{il}  A_{k} - A_{kl}  A_{i} + 2  A_{i }  A_{k}  A_{l}) \label{extra1} \\ 
& = & \! 4 \cdot  (A_{jkl} - A_{jk}  A_{l} - A_{jl}  A_{k} - A_{kl}  A_{j} + 2  A_{j }  A_{k}  A_{l})
  (A_i  A_j -A_{ij})(A_i  A_k-A_{ik})(A_i  A_l -A_{il}) . \nonumber
\end{eqnarray}
Thus conditions~(\ref{extra1}) are necessary for the realizability of $A_*$.
Note that the $2{\times}2{\times}2$ hyperdeterminantal relations  alone imply
a weaker version of~(\ref{extra1}), with both sides squared.

Also for $n\geq 4$ and for distinct $i$, $j$, $k$, $l$, we define
\begin{eqnarray*}
f_{ijkl}(A_*) & \eqbd  & A_\emptyset^3 A_{ijkl} - 3 \cdot A_i A_j A_k A_l - A_\emptyset^2 A_i A_{ijk} - A_\emptyset^2 A_j A_{ikl} - A_\emptyset^2 A_k A_{ijl} - A_\emptyset^2 A_{l} A_{ijk}  \\
&& + A_\emptyset \cdot (A_i A_j A_{kl} + A_i A_k A_{jl} + A_i A_l A_{jk} + A_j A_k A_{il} + A_j A_l A_{ik} + A_k A_l A_{ij}) \\
&& - (A_i A_j - A_\emptyset A_{ij})(A_k A_l -  A_\emptyset A_{kl}) - (A_i A_k -  A_\emptyset A_{ik})(A_j A_l -  A_\emptyset A_{jl})\\ && -
(A_i A_l -  A_\emptyset A_{il})(A_j A_k -  A_\emptyset A_{jk}) .
\end{eqnarray*}
This operation extracts the products of entries of $A$ that occur in the determinant
$A_{ijkl}$ indexed by permutations of order four. Namely, we find that
$$  f_{ijkl}(A_*)
\quad = \quad -2 \cdot (a_{ij} a_{il} a_{jk} a_{kl} + a_{ij} a_{ik} a_{jl} a_{kl} + a_{ik} a_{il} a_{jk} a_{jl} ).  $$
To rewrite these terms  differently, we use the polynomials
$$ g_{ijk}(A_*) \,\,\, \eqbd \,\,\,  A_\emptyset^2 A_{ijk} - A_\emptyset
 \cdot (A_i A_{jk} + A_j A_{ik} + A_k A_{ij}) + 2 \cdot A_i A_j A_k \;
\,\, \, = \,\,\,2\cdot a_{ij} a_{ik} a_{jk} , $$
and we observe that
\begin{eqnarray*}
g_{ikl}(A_*) g_{ijk}(A_*) &  = &  4 \cdot a^2_{ik} a_{ij} a_{il} a_{jk} a_{kl}, \\
g_{ikl}(A_*) g_{ijl}(A_*) &  = &  4 \cdot a^2_{il} a_{ij} a_{ik} a_{jl} a_{kl}, \\
g_{ijl}(A_*) g_{ijk}(A_*) &  = &  4 \cdot a^2_{ij} a_{ik} a_{il} a_{jk} a_{jl}.
\end{eqnarray*}
Using the relation $a^2_{ij} = A_i A_j -A_{ij}$, we thus obtain
\begin{eqnarray}
&& 2 f_{ijkl}(A_*) (A_i A_j - A_\emptyset A_{ij}) (A_i A_k -A_\emptyset A_{ik}) (A_i A_l -
 A_\emptyset A_{il}) \nonumber \\
&& \qquad  + g_{ikl}(A_*) g_{ijk}(A_*) (A_i A_j -A_\emptyset A_{ij}) (A_i A_l -A_\emptyset A_{il}) \nonumber \\
&& \qquad  + g_{ikl}(A_*) g_{ijl}(A_*)(A_i A_j -A_\emptyset A_{ij}) (A_i A_k -A_\emptyset A_{ik}) \nonumber \\
&& \qquad + g_{ijl}(A_*) g_{ijk}(A_*)  (A_i A_k -A_\emptyset A_{ik}) (A_i A_l -A_\emptyset A_{il})
\,\,\, =\,\,\,  0.  \label{extra2}
\end{eqnarray}
By the Schur complement argument from Section~\ref{sec_background}, we conclude
that the relations~(\ref{extra2}) continue to hold after the substitution
\begin{equation}
\begin{array}{l}
A_{\emptyset} \mapsto A_{I}, \,
A_{i} \mapsto A_{I \cup \{i\}},\,\,
A_{j} \mapsto A_{I \cup \{j\}},\,
A_{k} \mapsto A_{I \cup \{k\}}, \,
A_{l} \mapsto A_{I \cup \{l\}}, \\
A_{ij} \mapsto A_{I \cup \{i,j\}},\,
A_{ik} \mapsto A_{I \cup \{i,k\}}, \,
A_{il} \mapsto A_{I \cup \{i,l\}}, \\
A_{jk} \mapsto A_{I \cup \{j, k\}}, \,
A_{jl} \mapsto A_{I \cup \{j,l\}}, \,
A_{kl} \mapsto A_{I \cup \{k,l\}}, \, \\
A_{ijk} \mapsto A_{I \cup \{i,j,k\}},\,
A_{ijl} \mapsto A_{I \cup \{i,j,l\}},\,
A_{ikl} \mapsto A_{I \cup \{i,k,l\}}, \,
 A_{jkl} \mapsto A_{I \cup \{j,k,l\}}, \\
A_{ijkl} \mapsto A_{I \cup \{i,j,k,l\}}
\end{array}  \label{substitute2}
\end{equation}
for all subsets $I \subseteq [n]\backslash \{i,j,k,l \}$. We thus proved
the following addition to Theorem~\ref{first}.

\begin{lemma}  \label{second}
The principal minors of a symmetric matrix satisfy the conditions~(\ref{extra1})
and all conditions obtained from~(\ref{extra2}) via the  substitution~(\ref{substitute2}).
\end{lemma}

We are finally in a position to prove our first converse to Theorem~\ref{first}.
We assume a nondegeneracy condition to the effect that
a subset of the Hadamard-Fischer conditions
holds strictly. The condition   (\ref{basicHF}) is equivalent to the
condition that {\em weak sign symmetry} 
(see, e.g.,~\cite{gkktau}) holds strictly.

\begin{theorem}  \label{converse1}
Let $A_*$ be a real vector of length $2^n$ with $A_\emptyset=1$
that satisfies
\begin{equation}
A_{I\cap J} A_{I\cup J} < A_I A_J \quad  \hbox{\rm whenever} \;\; \# I\cap J = \# I - 1 = \# J  - 1  \label{basicHF}.
\end{equation}
There exists a real symmetric matrix $A$ with principal minors given by $A_*$
if and only if $A_*$ satisfies the
hyperdeterminantal relations, the relations~(\ref{extra1}) and all relations
obtained from (\ref{extra2}) using~(\ref{substitute2}).
\end{theorem}

\begin{proof}
Assuming~(\ref{basicHF}), (\ref{extra1}) and (\ref{extra2}), we build the
real symmetric matrix $A = (a_{ij}) $
as follows. The entries $A_I$ indexed by sets $I$ of size $1$ and $2$ determine all diagonal
entries $a_{ii}$ and the magnitudes of all off-diagonal entries. Note that~(\ref{basicHF})
implies that all off-diagonal entries of $A$ are non-zero. It remains to choose
the signs of off-diagonal entries correctly. Since the principal minors of $A$ do not change
under  diagonal similarity $A\mapsto DAD^{-1}$ where each diagonal entry of the matrix $D$
is $\pm 1$, we can fix all first row entries $a_{1j}$ with $j>1$ to be positive. Then the sign of each
entry $a_{2j}$ with $j>2$ is determined unambiguously by the values $A_{12j}$, the sign
of each entry $a_{3j}$ with $j>3$ is determined by the values $A_{13j}$, and so on.
In this fashion we prescribe all entries of the matrix $A$.

Using hyperdeterminantal relations and condition~(\ref{extra1}), we see that
this assignment is consistent with the values of all principal minors $A_{I}$
where  the index set $I$ has size at most $3$. Indeed, the hyperdeterminantal
relations guarantee that the absolute values of all off-diagonal entries are
consistent with the values of all principal minors of order at most $3$,
whereas conditions~(\ref{extra1}) guarantee that the signs are consistent
as well. 
The remaining entries of $A_*$,
i.e., the principal minors indexed by sets of size more than $3$, are determined
uniquely by the already specified entries of $A_*$. This follows from the conditions~(\ref{extra2}) and
from all conditions obtained from~(\ref{extra2}) via the substitution~(\ref{substitute2}),
since in every such condition  its largest principal minor occurs linearly with a nonzero
coefficient (which is a function of smaller minors) due to the inequalities~(\ref{basicHF}).
Thus the constructed matrix $A = (a_{ij})$
has the prescribed vector of principal minors $A_*$.   \end{proof}

\section{The Prime Ideal for $4 {\times} 4$-Matrices}   \label{sec_ideal}

From the point of view of algebraic geometry,
the following version of our problem is most natural:

\begin{problem}
Let $P_n$ be the prime ideal of all
homogeneous polynomial relations among
the principal minors of a symmetric $n {\times} n$-matrix.
Determine a finite set of generators for
the prime ideal $P_n$.
\end{problem}

The ideal $P_n$ lives in the  ring of polynomials
in the $2^n$ unknowns $A_I$ with rational coefficients.
For $n=3$, the ideal $I_3$ is principal, and its generator
is  the $2 {\times} 2 {\times} 2$-hyperdeterminant.
In Theorem~\ref{first} we identified the
subideal $H_n$ which is generated by all
hyperdeterminantal relations, one for each
$3$-dimensional face of the $n$-cube.
For instance, the ideal $H_4$ for the $4$-cube
is generated by the  eight homogeneous polynomials
of degree four  in $16$ unknowns listed in Example~\ref{eighthyper}.
It can be shown that $P_n$ is a minimal prime
of the ideal $H_n$ but in general we do not
know the minimal generators of $P_n$.

Two important features of both  ideals $P_n$ and $H_n$
is that they are invariant under the symmetry
group of the $n$-cube, and they are
 homogeneous with respect to the
$(n+1)$-dimensional multigrading induced by
the $n$-cube.  Both features
were used to simplify and organize our
 computations.
In this section we focus on the case $n=4$,
for which we establish the following result.

\begin{theorem} \label{idealgens1}
The homogeneous prime ideal $P_4$ is minimally generated by
twenty quartics in the
$16$ unknowns $A_I$. The corresponding irreducible
variety in $\P^{15}$ has codimension five and degree $96$.
\end{theorem}

This theorem was established with the aid of
computations using the computer algebra
packages {\tt Singular} and {\tt Macaulay 2}.
We worked in  the polynomial ring with the $5$-dimensional multigrading
\begin{eqnarray*} &
{\rm deg}(A_\emptyset) = (1,0,0,0,0),\,\,
{\rm deg}(A_1) = (1,1,0,0,0),\,\ldots, \,
{\rm deg}(A_4) = (1,0,0,0,1) \\ &
{\rm deg}(A_{12}) = (1,1,1,0,0),\ldots,\,\,
{\rm deg}(A_{234}) = (1,0,1,1,1), \,\,
{\rm deg}(A_{1234}) = (1,1,1,1,1).
\end{eqnarray*}
The twenty minimal generators of $P_4$ come in three
symmetry classes, with respect to the symmetry group
of the $4$-cube (i.e.~the Weyl group $B_4$ of order $384$).
 The three symmetry classes are as follows:

\vskip .1cm

\noindent {\bf Class 1:}
The eight $2 {\times} 2 {\times} 2$ {\em hyperdeterminants}
are listed in Example \ref{eighthyper}. Their multidegrees are
$$ \begin{matrix}
(4,2,2,2,0), &
(4,2,2,0,2), &
(4,2,0,2,2), &
(4,0,2,2,2),\\
(4,2,2,2,4), &
(4,2,2,4,2), &
(4,2,4,2,2), &
(4,4,2,2,2).
\end{matrix}
$$

\vskip .1cm

\noindent {\bf Class 2:}  There is a unique (up to scaling) minimal generator
of $P_4$ in each of the eight multidegrees
$$ \begin{matrix}
(4,2,2,2,1), &
(4,2,2,1,2), &
(4,2,1,2,2), &
(4,1,2,2,2),\\
(4,2,2,2,3), &
(4,2,2,3,2), &
(4,2,3,2,2), &
(4,3,2,2,2).
\end{matrix}
$$
The representative for multidegree
$(4,3,2,2,2)$ is the following quartic with $40$ terms:
\begin{eqnarray*}
&
A_{1234}^2 A_{1} A_{\emptyset}
-A_{1234} A_{123} A_{14} A_{\emptyset}
-A_{1234} A_{123} A_{1} A_{4}
-A_{1234} A_{124} A_{13} A_{\emptyset}
-A_{1234} A_{124} A_{1} A_{3} \\ &
-A_{1234} A_{134} A_{12} A_{\emptyset}
-A_{1234} A_{134} A_{1} A_{2}
-A_{1234} A_{12} A_{34} A_{1}
-A_{1234} A_{13} A_{24} A_{1}
-A_{1234} A_{14} A_{23} A_{1} \\ &
-A_{123} A_{124} A_{13} A_{4}
-A_{123} A_{124} A_{14} A_{3}
-A_{123} A_{134} A_{12} A_{4}
-A_{123} A_{134} A_{14} A_{2}
-A_{123} A_{234} A_{14} A_{1} \\ &
-A_{123} A_{12} A_{14} A_{34}
-A_{123} A_{13} A_{14} A_{24}
-A_{124} A_{134} A_{12} A_{3}
-A_{124} A_{134} A_{13} A_{2}
-A_{124} A_{234} A_{13} A_{1} \\ &
-A_{124} A_{12} A_{13} A_{34}
-A_{124} A_{13} A_{14} A_{23}
-A_{134} A_{234} A_{12} A_{1}
-A_{134} A_{12} A_{13} A_{24}
-A_{134} A_{12} A_{14} A_{23} \\ &
+2 A_{1234} A_{12} A_{13} A_{4}
+2 A_{1234} A_{12} A_{14} A_{3}
+2 A_{1234} A_{13} A_{14} A_{2}
+2 A_{123} A_{124} A_{134} A_{\emptyset}
+2 A_{123} A_{124} A_{34} A_{1} \\ &
+2 A_{123} A_{134} A_{24} A_{1}
+2 A_{124} A_{134} A_{23} A_{1}
+2 A_{234} A_{12} A_{13} A_{14}
+A_{1234} A_{234} A_{1}^2
+A_{123}^2 A_{14} A_{4} \\ &
+A_{123} A_{14}^2 A_{23}
+A_{124}^2 A_{13} A_{3}
+A_{124} A_{13}^2 A_{24}
+A_{134}^2 A_{12} A_{2}
+A_{134} A_{12}^2 A_{34}.
\end{eqnarray*}

\vskip .1cm

\noindent {\bf Class 3:}  The ideal $P_4$ contains
a four-dimensional space of minimal generators
in multidegree $(4,2,2,2,2)$. These generators
are not unique but we can choose a symmetric collection
of generators by taking the four quartics in the $B_4$-orbit of
the following polynomial with $36$ terms:
\begin{eqnarray*}
& 12 \cdot \bigl(
A_{\emptyset}^2 A_{1234}^2
+A_{4}^2 A_{123}^2
+A_{3}^2 A_{124}^2
+A_{34}^2 A_{12}^2
+A_{2}^2 A_{134}^2
+A_{24}^2 A_{13}^2
+A_{23}^2 A_{14}^2
+A_{234}^2 A_{1}^2 \bigr) \\ &
- A_{\emptyset} A_{34} A_{12} A_{1234}
- A_{\emptyset} A_{24} A_{13} A_{1234}
- A_{\emptyset} A_{23} A_{14} A_{1234}
- A_{\emptyset} A_{234} A_{1} A_{1234}
+ A_{\emptyset} A_{234} A_{14} A_{123} \\  &
+ A_{\emptyset} A_{234} A_{13} A_{124}
+ A_{\emptyset} A_{234} A_{134} A_{12}
- A_{4} A_{3} A_{124} A_{123}
- A_{4} A_{2} A_{134} A_{123}
+ A_{4} A_{23} A_{1} A_{1234} \\ &
- A_{4} A_{23} A_{14} A_{123}
+ A_{4} A_{23} A_{13} A_{124}
+ A_{4} A_{23} A_{134} A_{12}
- A_{4} A_{234} A_{1} A_{123}
- A_{3} A_{2} A_{134} A_{124} \\ &
+ A_{3} A_{24} A_{1} A_{1234}
+ A_{3} A_{24} A_{14} A_{123}
- A_{3} A_{24} A_{13} A_{124}
+ A_{3} A_{24} A_{134} A_{12}
- A_{3} A_{234} A_{1} A_{124} \\ &
+ A_{34} A_{2} A_{1} A_{1234}
+ A_{34} A_{2} A_{14} A_{123}
+ A_{34} A_{2} A_{13} A_{124}
- A_{34} A_{2} A_{134} A_{12} \\ &
- A_{34} A_{24} A_{13} A_{12}
- A_{34} A_{23} A_{14} A_{12}
- A_{2} A_{234} A_{1} A_{134}
- A_{24} A_{23} A_{14} A_{13}.
\end{eqnarray*}
It can be checked, using {\tt Macaulay 2} or {\tt Singular}, that
the higher degree polynomials
gotten from (\ref{extra1}) and (\ref{extra2})
by homogenization with $A_\emptyset$
are indeed polynomial linear combinations of these
$20$ quartics. We thus obtain the following strong converse to
Theorem \ref{first} in the case of $4 {\times} 4$-matrices:

\begin{corollary} \label{idealgens2}
A vector $A_* \in \C^{16}$ with $A_\emptyset = 1$ can be realized
as the principal minors of a symmetric matrix $A \in \C^{4 \times 4}$
if and only the above twenty quartics are zero at $A_*$. If the entries of the given vector
$A_*$ are real and satisfy (\ref{basicHF})
then the entries of the symmetric matrix $A$ are real numbers.
\end{corollary}

\begin{proof}
Consider the map which takes complex symmetric $4 {\times} 4$-matrices
to their vector of principal minors. The image of this map is closed by
Corollary \ref{closure}, and it hence equals the affine variety
in $\C^{15} = \{A_\emptyset = 1\}$ defined by the prime ideal $P_4$.
The statement for $\R$ is derived from Theorem \ref{converse1}.
\end{proof}

\section{Big Hyperdeterminants and Condensation Polynomials}   \label{sec_bighyperdet}

One ultimate goal is to generalize Theorem~\ref{idealgens1}
and Corollary~\ref{idealgens2} from $n=4$ to $n\geq 5$.
This section offers first steps in this direction, starting with the
general hyperdeterminantal constraints on $A_*$.

We recall from~\cite[Chap.14]{GKZbook} that
the {\em hyperdeterminant} of a tensor
$\,\A = (\a_{i_1, \ldots,i_r} )\,$
of format $\,2 {\times} 2 {\times} \cdots {\times} 2\,$
is defined as follows. Consider the multilinear form
$f$ defined by the tensor $\A$:
\begin{equation}
 f(x) \,\,\eqbd \,\, f(x^{(1)}, x^{(2)} \ldots, x^{(r)}) \,\,\,\eqbd\,\,\, \sum_{i_1=0}^{1}  \sum_{i_2=0}^{1}
  \cdots \sum_{ i_r=0}^{1} \a_{i_1,i_2, \ldots,i_r} \cdot
x^{(1)}_{i_1} x^{(2)}_{i_2} \cdots x^{(r)}_{i_r}.  \label{main_form}
\end{equation}
The hyperdeterminant $\det(\A)$ is the
unique (up to scaling) irreducible polynomial in the entries of $\A$ that
characterizes the degeneracy of the form $f$, i.e., $\det(\A)=0$ if and only if the equations
\begin{equation}
f(x)\,\,= \,\, {\partial f(x) \over \partial x_i^{(j)}}
\,\, = \,\, 0 \quad \hbox{\rm for all}\;\; i, \, j  
\label{degener}
\end{equation}
have a solution $x=(x^{(1)},x^{(2)}, \ldots, x^{(r)})$ where each $x^{(j)}$
is a non-zero complex vector in $\C^2$.

As before we identify our proposed vector of principal minors
 $A_* \in \R^{2^n}$ with the corresponding $2 {\times} 2 {\times} \cdots {\times} 2$-tensor.
The following result generalizes Theorem~\ref{first} and gives an alternative proof.

\begin{theorem} Let $A = (a_{ij})$ be a symmetric $n{\times}n$ matrix. Then the tensor $A_*$ of
all principal minors of $A$ is a common zero of all the hyperdeterminants of formats
from $2 {\times} 2 {\times} 2$ up to
$\two_n$. \end{theorem}

\begin{proof} It suffices to prove that the highest-dimensional hyperdeterminant
vanishes, using Schur complements and
 induction. Let $f$ be the form~(\ref{main_form})
corresponding to the tensor $A_*$. Take
\begin{equation}
  x^{(1)} \,\eqbd\, ( a_{12}a_{13}-a_{11}a_{23} \, , \, a_{23} ), \quad
 x^{(2)} \,\eqbd \, (a_{12}a_{23}-a_{13}a_{22}\, , \, a_{13}),
 \quad x^{(3)} \,\eqbd \,(a_{13}a_{23}-a_{12}a_{33}\, , \, a_{12}).  \label{x_chosen}
\end{equation}
 and take the remaining $x^{(j)}$ to be $(1,0)$.
With this choice of $x$, the conditions~(\ref{degener}) are satisfied.
\end{proof}

We next introduce the {\em condensation polynomial} $C_n$
which expresses the determinant $A_{123 \cdots n}$
as an algebraic function of the principal minors
$A_i$ and $A_{ij}$ of size at most two.
For instance, for $n=3$ the determinant $A_{123}$ is
an algebraic function of degree two in $\{A_1,A_2,A_3,A_{12}, A_{13}, A_{23}\}$,
and thus $C_3$ coincides with  the $2{\times} 2 {\times} 2$-hyperdeterminant.
In general, the condensation polynomial $C_n$ is defined as the
unique irreducible (and monic in $A_{123\cdots n }$) generator of the
principal elimination ideal
$$ \langle C_n \rangle \quad = \quad
(P_n + \langle A_\emptyset -1 \rangle ) \, \cap \,
\Q[\,A_1, A_2,\ldots, A_n, \,A_{12} ,A_{13} \ldots, A_{n-1,n} \, ,
\, A_{123\cdots n} \,] . $$
Using the sign-swapping argument in the proof of
Theorem~\ref{converse1}, we can show that
$C_n$ is a polynomial of degree $2^{\binom{n-1}{2}}$
in $A_{123\cdots n}$. The total degree of $C_n$
is bounded above by $\,n \cdot 2^{\binom{n-1}{2}}$.
The following derivation proves these assertions for $n=4$,
and it illustrates the construction of $C_n$ in general.

\begin{example} \label{condensation} \rm
The condensation polynomial $C_4$ is an irreducible
polynomial of degree $23$. It is the sum of $ 12380$ monomials in
$\, \Q[ A_1, A_2, A_3,A_4,
A_{12},A_{13}, A_{14},A_{23},A_{24},A_{34},A_{1234}]$.
To compute $C_4$, we first consider  the following
polynomial which expresses the symmetric
$4 {\times} 4$-determinant
\begin{equation}
\label{4x4det}
A_{1234} \, - \, {\rm det} \begin{pmatrix} A_1 & a_{12} & a_{13} & a_{14} \\
                                            a_{12} & A_2 & a_{23} & a_{24} \\
                                            a_{13} & a_{23} & A_3 & a_{34} \\
                                             a_{14} & a_{24} & a_{34} & A_4
\end{pmatrix} .
\end{equation}
To get rid of the square roots
$\,a_{ij} = \sqrt{A_i A_j - A_{ij} } $,
we swap the sign on the $a_{ij}$ with $1 < i < j$
 in all eight possible ways.
The orbit of (\ref{4x4det}) under these sign swaps
consists of eight distinct variants of (\ref{4x4det}). The product of these
eight expressions equals $C_4$. Interestingly, the
degree drops to $23$.
 \qed
\end{example}

We could in fact use the condensation polynomials $C_n$
to replace the hypothesis~(\ref{basicHF}) from the
statement of Theorem~\ref{converse1}, thus
providing a stronger converse to Theorem~\ref{first}.
However, this is not particularly useful for a practical test,
since the condensation polynomials $C_n$ are too big to compute
explicitly for $n \geq 5$. What is desired instead are
explicit  generators of the prime ideal $P_n$.

\section{Invariance and the Quartic Generation Conjecture} \label{sec_landsberg}

This section was written after the first version of this paper
had been submitted for publication. It is based on discussions
with J.M.~Landsberg, who suggested Theorem \ref{invariance} to us in May 2006.

Let $\R[A_\bullet]$ denote the polynomial ring in $2^n$ unknowns
$A_I$ where $I$ runs over all subsets of $\{1,2,\ldots,n\}$.
We are interested in the  prime ideal $P_n$ of all homogeneous
polynomials in $\R[A_\bullet]$ which are algebraic relations among
the $2^n$ principal minors of a generic symmetric $n \times n$-matrix.
Clearly, $P_n$ is invariant under the action of the symmetric group
$S_n$ on the polynomial ring $\R[A_\bullet]$. In what follows we show
that $P_n$ is invariant under a natural Lie group action
on $\R[A_\bullet]$ as well.  Let $SL_2(\R) $ denote group of
real $2 \times 2$-matrices with determinant $1$. The $n$-fold product of
this group,
$$ G \quad := \quad SL_2 (\R) \times SL_2(\R) \times \cdots \times SL_2 (\R), $$
acts naturally on the $n$-fold tensor product
$$ \R^{2^n} \quad := \quad  \R^2 \otimes \R^2 \otimes \cdots \otimes \R^2. $$
Since $\R[A_\bullet]$ is the ring of polynomial functions of $\R^{2^n}$,
we get an action of $G$ on $\R[A_\bullet]$.

\begin{theorem} \label{invariance}
The homogeneous prime ideal $P_n$ is invariant under the group $G$.
\end{theorem}

\begin{proof}
Consider an $n \times 2n$-matrix $\,(B \,\,C)\,$ where
$B$ and $C$ are generic $n \times n$-matrices subject
to the constraint that $B^{-1} C$ is symmetric.
We identify each principal minor of the symmetric matrix $B^{-1} C$
with an $n \times n$-minor of the matrix $\,(B \,\, C )$.
The $2^n$ maximal minors of $\,(B \,\, C )\,$
which come from principal minors of $B^{-1}C$ are
precisely those whose column index sets $J \subset \{1,2,\ldots,2n\}$ satisfy
$\# (J \,\cap \,\{i,i+n\}) = 1\,$ for $i=1,2,\ldots,n$.
We call these the {\em special maximal minors} of $(B \,\, C)$.
The prime ideal $P_n$ consists of the algebraic relations among
the $2^n$ special maximal minors.

We need to show that the parametric variety corresponding
to the prime ideal $P_n$ is invariant under the group $G$.
In order to do this, we consider the representation of the group $G$
by $2n \times 2n$-matrices of the form
$$ g \,\,\, = \,\,\, \begin{pmatrix} D_1 & D_2 \\ D_3 & D_4 \end{pmatrix}$$
where $D_1,D_2,D_3,D_4$ are diagonal $n{\times}n$-matrices which satisfy the identity
$$ \qquad \qquad D_1 \cdot D_4 \, - \,D_2 \cdot D_3 \quad = \quad {\bf 1} \qquad \qquad \hbox{
(the $n{\times}n$-identity matrix}) $$
Here, the $i$-th factor $SL_2(\R)$ in the $n$-fold product $G$ corresponds to
the $2{\times}2$-matrix formed by the entries in position $(i,i)$
of the diagonal matrices $D_1,D_2,D_3,D_4$.
Let $g$ be any element of $G$ represented by a $2n{\times}2n$-matrix.
Then the vector in $\R^{2^n}$ of  special maximal minors of $\,(B \,\, C) \cdot g \,$
is precisely the result of applying $g$ to the vector
in $\R^{2^n}$ of special maximal minors of $\,(B \,\, C)$.

It now suffices to show is the following matrix-theoretic statement:
if $B,C,\tilde B, \tilde C$ are $n \times n$-matrices such that
 $\,(\tilde B\,\, \tilde C) \,= \, (B \,\, C) \cdot g \,$ for some $g \in G$, and if
 $B^{-1} C $ is symmetric, then
$ \tilde B^{-1} \tilde C$ is also symmetric. The following lemma
proves this statement. \end{proof}

\begin{lemma}
Let $B$ and $ C$ be invertible $n \times n$-matrix
such that $B^{-1} C$ is symmetric,
and let $D_1,D_2,D_3,D_4$ be diagonal $n \times n$-matrices
satisfying $\, D_1 \cdot D_4 \, - \,D_2 \cdot D_3 \, = \, {\bf 1} $.
Then the matrix $\,(B D_1 + C D_2)^{-1} (B D_3 + C D_4)\,$ is also symmetric.
\end{lemma}

\begin{proof}
The assumption that $B^{-1} C$ is symmetric is equivalent to
the following identity:
\begin{equation}
\label{symmassum}
 B \cdot C^T  \,\, - \,\, C \cdot B^T \quad = \quad {\bf 0}.
\end{equation}
Similarly, the desired conclusion states that the following difference is the zero matrix ${\bf 0}$:
$$ (B D_1 + C D_2) \cdot (B D_3 + C D_4)^T \,- \,
(B D_3 + C D_4) \cdot (B D_1 + C D_2)^T $$
Multiplying out, cancelling common terms, and  using both
$D_i^T = D_i$ and the identity (\ref{symmassum}), we simplify the
above matrix expression as follows:
\begin{eqnarray*}
& C D_2 D_3 B^T + B D_1 D_4 C^T - B D_3 D_2 C^T - C D_4 D_1 B^T  \\
= &
B \cdot (D_1 D_4 - D_2 D_3 ) \cdot C^T \,-\, C \cdot  (D_1 D_4 - D_2 D_3) \cdot B^T \\
= & B {\bf 1} C^T - C {\bf 1} B^T \quad = \quad B C^T - C B^T \quad = \quad {\bf 0}.
\end{eqnarray*}
This proves the lemma and hence the theorem.
\end{proof}

We define the  {\em hyperdeterminantal module}
to the $G$-orbit of the $2 {\times} 2 {\times} 2$-hyperdeterminant
under the action of the group $G$ and the symmetric group $S_n$.
This orbit is a subspace of the finite-dimensional vector space
$\,\R[A_\bullet]_4\,$ of quartic polynomials in the $2^n$ unknowns $A_I$.
Using a representation theoretic argument, it can be shown that
the vector space dimension of the hyperdeterminantal module equals
$$ \binom{ 2^{n-3} + 3}{ 4 } \cdot \binom{n}{3} $$
This number is one for $n=3$, and it is $20$ for $n=4$.
We propose the following natural conjecture.

\begin{conjecture}
The prime ideal $P_n$ is generated by the hyperdeterminantal module.
\end{conjecture}

For $n=3$, the prime ideal is principal and generated by the
$2 {\times} 2 {\times} 2$-hyperdeterminant.
For $n = 4$, this conjecture is established by our computations in Section 4.

\section*{Acknowledgements}

We thank David Wagner and Yuval Peres for information on the probabilistic aspects
of the principal minor assignment problem.
Theorem~\ref{invariance} was suggested to us by J.M.~Landsberg,
in response to the first version of this paper.
We thank him for his permission to include this result here.
Bernd Sturmfels was supported by the National Science Foundation
(DMS-0456960).

\bibliographystyle{plain}

\def\cprime{$'$} \def\cprime{$'$} \def\cprime{$'$}

\end{document}